\documentclass[12pt]{article}

\usepackage{amsmath,amssymb,
amsthm,graphicx,extsizes}
\title{Classification of sesquilinear
forms with the first argument on
a subspace or a factor space\footnotetext{This is the authors' version of a work that was published in Linear Algebra Appl. 424 (2007) 282--303.}}
\author{Vyacheslav Futorny%
\thanks{The author was partially supported by CNPq,
processo 307812/2004-9.}
\\Department of Mathematics, University of
S\~{a}o Paulo\\ S\~{a}o Paulo,
Brazil, futorny@ime.usp.br
 \and
Vladimir V. Sergeichuk%
\thanks{Corresponding author. The research was
done while this author was
visiting the University of
S\~{a}o Paulo supported by
FAPESP, processo 05/59407-6.}
\\ Institute of
Mathematics, Tereshchenkivska
3\\ Kiev, Ukraine,
sergeich@imath.kiev.ua}
\date{\it Dedicated to R. A. Horn
on the occasion of his 65th
birthday}

\begin{document}

 \maketitle
\begin{abstract}
Let $V$ be a vector space over a
field or skew field $\mathbb F$,
and let $U$ be its subspace. We
study the canonical form problem
for bilinear or sesquilinear
forms
\[
U\times V\rightarrow {\mathbb
F},
  \qquad
(V/U)\times V\rightarrow
{\mathbb F}
\]
and linear mappings $
U\rightarrow V,$ $V\rightarrow
U,$  $V/U\rightarrow V,$
$V\rightarrow V/U.$ We solve it
over $\mathbb F=\mathbb C$ and
reduce it over all $\mathbb F$
to the canonical form problem
for ordinary linear mappings
$W\to W$ and bilinear or
sesquilinear forms $W\times
W\rightarrow {\mathbb F}$.
Moreover, we give an algorithm
that realizes this reduction.
The algorithm uses only unitary
transformations if $\mathbb
F=\mathbb C$, which improves its
numerical stability. For linear
mapping this algorithm can be
derived from the algorithm by L.
A. Nazarova, A. V. Roiter, V. V.
Sergeichuk, and V. M. Bondarenko
[\emph{J. Soviet Math.} 3 (no.
5) (1975) 636--654].


\end{abstract}

\newcommand{\ddd}{
\text{\begin{picture}(12,8)
\put(-2,-4){$\cdot$}
\put(3,0){$\cdot$}
\put(8,4){$\cdot$}
\end{picture}}}

\newcommand{\dddd}{
\text{\begin{picture}(12,8)
\put(-2,-3){$\cdot$}
\put(3,0){$\cdot$}
\put(8,3){$\cdot$}
\end{picture}}}

\newcommand{\sdotsss}%
{\text{\raisebox{-2.2pt}{$\cdot\,$}%
 \raisebox{1.7pt}{$\cdot$}%
\raisebox{5.6pt}{$\,\cdot$}}}

\renewcommand{\le}{\leqslant}
\renewcommand{\ge}{\geqslant}

\newcommand{\sdotss}%
{\text{\raisebox{-1.8pt}{$\cdot\,$}%
 \raisebox{1.5pt}{$\cdot$}%
\raisebox{4.8pt}{$\,\cdot$}}}

\newcommand{\diag}{\,\diagdown\,}
\newcommand{\rank}{\mathop{\rm rank}\nolimits}

\newtheorem{theorem}{Theorem}
\newtheorem{lemma}[theorem]{Lemma}

\theoremstyle{definition}
\newtheorem{definition}[theorem]{Definition}

\section{Introduction} \label{s-res}

In this paper, we give canonical
matrices of bilinear or
sesquilinear forms
\[
U\times V\rightarrow {\mathbb
C},
  \qquad
(V/U)\times V\rightarrow
{\mathbb C},
\]
where $V$ is a complex vector
space and $U$ is its subspace.

We use the following canonical
matrices of bilinear or
sesquilinear forms on a complex
vector space given in
\cite{hor-ser_transpose} (see
also \cite{hor-ser1_singul,
hor-ser1_canon}). Two square
complex matrices $A$ and $B$ are
said to be \emph{congruent} or
*\!\emph{congruent} if there is
a nonsingular $S$ such that
$S^{T}AS=B$ or, respectively,
$S^{*}AS=B$, where
$S^{*}:=\bar{S}^{T}$ denotes the
complex conjugate transpose of
$S$. Define the $n$-by-$n$
matrices
\begin{equation*}\label{1aa}
\Gamma_n = \begin{bmatrix}
0&&&&\dddd\\ &&&-1& \dddd\\
&&1&1&\\ &-1&-1& &\\ 1&1&&&0
\end{bmatrix},\qquad
\Delta_n=\begin{bmatrix} 0&&&1
\\
&&\ddd&i\\
&1&\ddd&\\
1&i&&0
\end{bmatrix},
\end{equation*}
\begin{equation*}\label{bip1}
J_n(\lambda) =
\begin{bmatrix}
  \lambda&1&&0\\
  &\lambda&\ddots&\\
  &&\ddots&1\\
  0&&&\lambda
\end{bmatrix}.
\end{equation*}

\begin{theorem}
[{\cite[p.\,351]{hor-ser_transpose}}]
\label{bilin} {\rm(a)} Every
square complex matrix is
congruent to a direct sum,
determined uniquely up to
permutation of summands, of
matrices of the form
\begin{equation*}\label{eqqz}
J_n(0),
 \quad
\Gamma_n,
 \quad
\begin{bmatrix}0&I_n\\ J_n(\lambda)
&0
\end{bmatrix},
\end{equation*}
in which $\lambda \ne 0$,
$\lambda\ne (-1)^{n+1}$, and
$\lambda$ is determined up to
replacement by $\lambda^{-1}$.

{\rm(b)} Every square complex
matrix is {\rm *}\!congruent to
a direct sum, determined
uniquely up to permutation of
summands, of matrices of the
form
\begin{equation*}\label{eqq}
J_n(0),
 \quad
\lambda\Gamma_n,
 \quad
\begin{bmatrix}0&I_n\\ J_n(\mu)
&0
\end{bmatrix},
\end{equation*}
in which $|\lambda|=1$ and\/
$|\mu|>1$. Alternatively, one
may use the symmetric matrix
$\Delta_n$ instead of\/
$\Gamma_n$. \hfill$\square$
\end{theorem}

A canonical form of a square
matrix for
congruence/*congruence over any
field $\mathbb F$ of
characteristic different from 2
was given in
\cite{ser_izvestiya} up to
classification of Hermitian
forms over finite extensions of
$\mathbb F$.

Let us formulate the main
result. For generality, we will
consider matrices over any field
or skew field $\mathbb F$ with
involution $\alpha\mapsto
\bar{\alpha}$, that is, a
bijection on $\mathbb F$ such
that
\begin{equation*}\label{1ig}
\overline{\alpha+\beta}=
\bar{\alpha}+\bar{\beta},
  \qquad
\overline{\alpha\beta}
=\bar{\beta} \bar{\alpha},
  \qquad
\bar{\bar{\alpha}}=\alpha
\end{equation*}
for all $\alpha,\beta\in \mathbb
F$.

We denote the $m$-by-$n$ zero
matrix by $0_{mn}$, or by $0_m$
if $m=n$. It is agreed that
there exists exactly one matrix
of size $n\times 0$ and there
exists exactly one matrix of
size $0\times n$ for every
nonnegative integer $n$; they
represent the linear mappings
$0\to {\mathbb F}^n$ and
${\mathbb F}^n\to 0$ and are
considered as the zero matrices
$0_{n0}$ and $0_{0n}$. For every
$p\times q$ matrix $M_{pq}$ we
have
\[
M_{pq}\oplus
0_{m0}=\begin{bmatrix}
  M_{pq} & 0 \\
  0 &0_{m0}
\end{bmatrix}=\begin{bmatrix}
  M_{pq}& 0_{p0} \\
  0_{mq}& 0_{m0}
\end{bmatrix}=\begin{bmatrix}
M_{pq} \\ 0_{mq}
\end{bmatrix}
\]
and
\[ M_{pq}\oplus
0_{0n}=\begin{bmatrix}
  M_{pq} & 0 \\
  0 & 0_{0n}
\end{bmatrix}=\begin{bmatrix}
  M_{pq}& 0_{pn} \\
  0_{0q}& 0_{0n}
\end{bmatrix}=\begin{bmatrix}
   M_{pq} & 0_{pn}
\end{bmatrix}.
\]
In particular, \[0_{p0}\oplus
0_{0q}=0_{pq}.\]

For each matrix $A=[a_{ij}]$
over $\mathbb F$, we define its
\emph{conjugate transpose}
\[
A^*=\overline{A}^{\mathrm
T}=[\bar{a}_{ji}].
\]
If $S^*AS=B$ for some
nonsingular matrix $S$, then $A$
and $B$ are said to be
*{\it\!congruent} (or
\emph{congruent} if $\mathbb F$
is a field and the involution on
$\mathbb F$ is the identity---in
what follows we consider
congruence as a special case of
*congruence).

A \emph{sesquilinear form} on
right vector spaces $U$ and $V$
over $\mathbb F$ is a map
\begin{equation*}\label{0.1}
{\cal G}\colon U\times
V\rightarrow \mathbb F
\end{equation*}
satisfying
\begin{align*}\label{0.2}
{\cal G}(u\alpha+u'\beta,v)
&=\bar{\alpha} {\cal G}(u,v)+
\bar{\beta} {\cal G}(u',v),\\
{\cal G}(u,v\alpha+v'\beta) &=
{\cal G}(u,v)\alpha +{\cal
G}(u,v')\beta
\end{align*}
for all $u,u'\in U,\ v,v'\in V$,
and $\alpha,\beta\in \mathbb F$.
If $\mathbb F$ is a field and
the involution on $\mathbb F$ is
the identity, then a
sesquilinear form becomes
bilinear---we consider bilinear
forms as a special case of
sesquilinear forms.

If $e_1,\dots,e_m$ and
$f_1,\dots,f_n$ are bases of $U$
and $V$, then
\begin{equation}\label{0.3}
G_{ef}=[\alpha_{ij}],\qquad
\alpha_{ij}:={\cal G}(e_i,f_j),
\end{equation}
is the \emph{matrix of ${\cal
G}$} in these bases. Its matrix
in other bases $e'_1,\dots,e'_m$
and $f'_1,\dots,f'_n$ can be
found by the formula
\begin{equation}\label{0.5}
G_{e'f'}=S^*G_{ef}R,
\end{equation}
where $S$ and $R$ are the change
of basis matrices.

For every $u\in U$ and $v\in V$,
\begin{equation*}\label{0.4}
{\cal G}(u,v)
=[u]_e^*\,G_{ef}\,[v]_f,
\end{equation*}
where $[u]_e$ and $[v]_f$ are
the coordinate column-vectors of
$u$ and $v$.
\medskip

In this paper, we study
sesquilinear forms
\begin{equation}\label{0.0}
U\times V\rightarrow {\mathbb
F}, \qquad (V/U)\times
V\rightarrow {\mathbb F},
\end{equation}
in which $U$ is a subspace of
$V$, so we always consider their
matrices \eqref{0.3} in those
bases of $U$ and $V$ that are
concordant as follows.

\begin{definition}\label{def0}
Let ${\cal G}$ be one of
sesquilinear forms \eqref{0.0},
in which $V$ is a right space
over $\mathbb F$, and $U$ is its
subspace. Choose a basis
$e_1,\dots,e_n$ of\/ $V$ such
that
\begin{equation}\label{0.21}
\begin{cases}
 \text{$e_1,\dots,e_m$ is a basis of $U$}
 &\text{if
${\cal G}\colon U\times V\to
{\mathbb F}$,}
   \\
\text{$e_{m+1},\dots,e_n$ is a
basis of $U$}
 &\text{if ${\cal G}\colon (V/U)\times V\to
{\mathbb F}$.}
\end{cases}
\end{equation}
By the \emph{matrix of ${\cal
G}$ in the basis
$e_1,\dots,e_n$}, we mean the
block matrix
\begin{equation}\label{0.22}
[A|B]=\left.\left[
\begin{matrix}
\alpha_{11}&\dots&\alpha_{1m}\\
\vdots&\ddots&\vdots\\
\alpha_{m1}&\dots&\alpha_{mm}
\end{matrix}
\right|
\begin{matrix}
\alpha_{1,m+1}&\dots&\alpha_{1n}\\
\vdots&&\vdots\\
\alpha_{m,m+1}&\dots&\alpha_{mn}
\end{matrix}
\right],
\end{equation}
in which
\begin{equation*}\label{0.23}
\alpha_{ij}=
  \begin{cases}
{\cal G} (e_i,e_j) &
    \text{if ${\cal
G}\colon U\times V\to {\mathbb
F}$},
        \\
{\cal G} (e_i+U,e_j) &
    \text{if ${\cal
G}\colon (V/U)\times V\to
{\mathbb F}$}.
  \end{cases}
\end{equation*}
\end{definition}

By the \emph{block-direct sum}
of block matrices $[A_1|B_1]$
and $[A_2|B_2]$, we mean the
block matrix
\begin{equation*}%\label{0.15}
[A_1|B_1]\uplus[A_2|B_2]:=
\left.\left[
\begin{matrix}
A_1&0\\0&A_2
\end{matrix}
\right|
\begin{matrix}
B_1&0\\0&B_2\end{matrix}
\right].
\end{equation*}

In Section \ref {s-pr} we will
prove the following theorem (a
stronger statement was proved in
\cite[Theorem
1]{hor-ser1_singul} in the case
$U=V$).

\begin{theorem}\label{t0.01}
Let $\mathbb F$ be a field or
skew field with involution
$($possibly, the identity if $F$
is a field$)$, $V$ be a right
vector space over $\mathbb F$,
and $U$ be its subspace. Let
${\cal G}$ be one of
sesquilinear forms
\begin{equation}\label{0.20}
U\times V\rightarrow {\mathbb
F}, \qquad (V/U)\times
V\rightarrow {\mathbb F}.
\end{equation}

{\rm(a)} There exists a basis
$e_1,\dots,e_n$ of\/  $V$
satisfying \eqref{0.21}, in
which the matrix \eqref{0.22} of
${\cal G}$ is a block-direct sum
of a $p$-by-$p$ matrix
\begin{equation}\label{0.10}
[K|0_{p0}],\qquad \text{$K$ is
nonsingular,}
\end{equation}
and matrices of the form
\begin{equation}\label{0.11}
[J_q(0)|0_{q0}]\ \ (q\ge 1),
 \qquad
[J_q(0)|E_q]\ \ (q\ge 0),
\end{equation}
in which
\begin{equation}\label{bip1s}
E_q:=\begin{bmatrix} \,0\,\\
\vdots\\ 0\\ 1
\end{bmatrix}\ \text{if }q\ge 1,
     \qquad
    E_0:=0_{01}
\end{equation}
$($the summands \eqref{0.10} or
\eqref{0.11} may be absent$)$.
The block $K$ is determined by
${\cal G}$ uniquely up to {\rm
*\!}congruence, and the summands
of the form \eqref{0.11} are
determined by ${\cal G}$
uniquely up to permutation.

 {\rm(b)}
If\/ $\mathbb F=\mathbb C$, then
one can replace in this direct
sum the summand \eqref{0.10} by
\begin{equation*}\label{mjs}
[K_1|0_{p_10}]\uplus\dots\uplus
[K_s|0_{p_s0}],
\end{equation*}
where $K_1\oplus\dots\oplus K_s$
is the canonical form of $K$
defined in Theorem \ref{bilin}
and each $K_i$ is
$p_i$-by-$p_i$. The obtained
block-direct sum is determined
by ${\cal G}$ uniquely up to
permutation of summands, and so
it is a canonical matrix of the
sesquilinear $($in particular,
bilinear\emph{)} form ${\cal
G}$.
\end{theorem}

Let us formulate an analogous
statement for matrices of linear
mappings.

\begin{definition}\label{defm}
Let $\mathbb F$ be a field or
skew field, $V$ be a right
vector space over $\mathbb F$,
and $U$ be its subspace. Let
${\cal A}$ be one of linear
mappings
\[
U\rightarrow V,\qquad
V\rightarrow U,\qquad
V/U\rightarrow V,\qquad
V\rightarrow V/U.
\]
Choose a basis $e_1,\dots,e_n$
of\/ $V$ such that
\begin{equation}\label{m.2}
\begin{cases}
 \text{$e_1,\dots,e_m$ is a basis of $U$,}
 &\text{if
$U\to V$ or $V\to U$,}
   \\
\text{$e_{m+1},\dots,e_n$ is a
basis of $U$,} &\text{if $V/U\to
V$ or
 $V\to V/U$.}
\end{cases}
\end{equation}
By the \emph{matrix $A_e$ of
${\cal A}$ in the basis
$e_1,\dots,e_n$}, we mean its
matrix in the bases
\[
\begin{cases}
 \text{$e_1,\dots,e_m$ of $U$,}
  &\text{if
$U\to V$ or $V\to U$,}
   \\
\text{$e_1+U,\dots,e_m+U$ of
$V/U$,}
  &\text{if $V/U\to V$ or
 $V\to V/U$,}
\end{cases}
\]
and $e_1,\dots,e_n$ of\/ $V$. We
divide $A_e$ into two blocks
\begin{equation}\label{m.3}
A_e= \begin{cases}
    \left[ \begin{array}{cc}
  A \\ \hline B
 \end{array}\right], & \text{if
$U\to V$ or $V/U\to V$,} \\[5mm]
   [A| B],
 & \text{if $V\to U$ or
 $V\to V/U$,}
  \end{cases}
\end{equation}
where $A$ is $m$-by-$m$.
\end{definition}

The following theorem will be
proved in Section  \ref{s-pr}.

\begin{theorem}\label{t.m}
Let $\mathbb F$ be a field or
skew field, $V$ be a right
vector space over $\mathbb F$,
and $U$ be its subspace. Let
${\cal A}$ be one of linear
mappings
\begin{equation}\label{m.4}
U\rightarrow V,\qquad
V\rightarrow U,\qquad
V/U\rightarrow V,\qquad
V\rightarrow V/U.
\end{equation}

{\rm(a)} There exists a basis
$e_1,\dots,e_n$ of\/ $V$
satisfying \eqref{m.2}, in which
for the matrix $A_e$ of ${\cal
A}$ we have:
\[
\begin{cases}
 A_e^T,
&\text{if $U\to V$ or $V/U\to
V$,}
   \\
A_e, &\text{if $V\to U$ or
 $V\to V/U$}
\end{cases}
\]
is a block-direct sum of a
$p$-by-$p$ matrix
\begin{equation}\label{m.5}
[K|0_{p0}],\qquad \text{$K$ is
nonsingular,}
\end{equation}
and matrices of the form
\begin{equation}\label{m.6}
[J_q(0)|0_{q0}],
 \qquad
[J_q(0)|E_q],
\end{equation}
where $E_q$ was defined in
\eqref{bip1s} $($the summands
\eqref{m.5} or \eqref{m.6} may
be absent$)$. The block $K$ is
determined by ${\cal A}$
uniquely up to similarity, and
the summands of the form
\eqref{m.6} are determined by
${\cal A}$ uniquely up to
permutation.

 {\rm(b)}
If\/ $\mathbb F=\mathbb C$, then
one can replace the summand
\eqref{m.5} by a block-direct
sum of square matrices of the
form
\begin{equation*}\label{m.7}
[J_q(\lambda)|0_{q0}].
\end{equation*}
The obtained matrix is
determined by ${\cal A}$
uniquely up to permutation of
summands, and so it is a
canonical matrix of the linear
mapping ${\cal A}$.
\end{theorem}

We do not rate Theorem \ref{t.m}
as new; it is readily available
from the canonical form problem
solved in \cite[\S\,2]{naz}. We
include it in our paper since
the singular indecomposable
summands of the canonical forms
in Theorems \ref{t0.01} and
\ref{t.m} coincide, and our
proofs of Theorems \ref{t0.01}
and \ref{t.m} are similar and
are based on
\emph{regularization algorithms}
that decompose the matrix of
each form \eqref{0.0} and each
mapping \eqref{m.4} into a
block-direct sum of
\begin{itemize}
  \item
 its \emph{regular part} $[K|0_{p0}]$
with nonsingular $K$ (see
\eqref{0.10} and \eqref{m.5}),
which is determined by
\eqref{0.0} or \eqref{m.4} up to
*congruence or similarity, and
of
  \item
 its \emph{singular summands} of the
form $[J_q(0)|0_{q0}]$ and
$[J_q(0)|E_q]$ (see \eqref{0.11}
and \eqref{m.6}), which are
determined uniquely.
\end{itemize}
If $\mathbb F=\mathbb C$, then
these algorithms can use only
unitary transformations, which
improves their numerical
stability. These algorithms
extend the regularization
algorithm \cite{hor-ser1_singul}
for a bilinear/sesquilinear
form, which decomposes its
matrix into a direct sum of a
nonsingular matrix and several
singular Jordan blocks. An
analogous regularization
algorithm was given by Van
Dooren \cite{doo} for matrix
pencils and was extended to
matrices of cycles of linear
mappings in \cite{ser}.

The canonical form problems for
matrices of forms \eqref{0.0}
and mappings \eqref{m.4} are
special cases of the canonical
form problem for block matrices,
whose form resembles
\[\includegraphics[width=1.5in,height=1.2in]
 %{bangle.jpg}\]
 %\[\includegraphics[bb=0 0 100 %80]
 {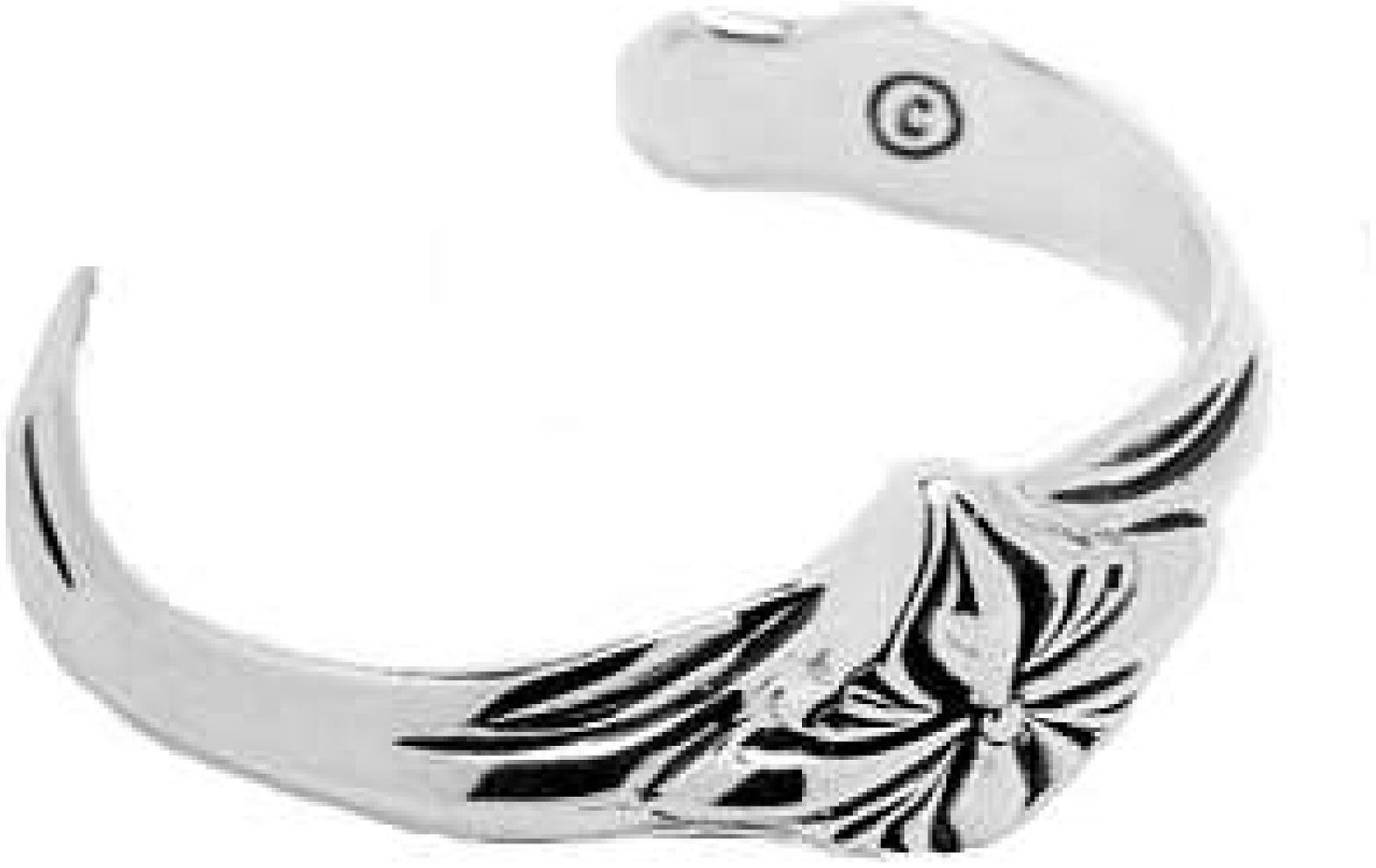}\]

\begin{definition} \label{def1}
By a \emph{bangle} over $\mathbb
F$ we mean a matrix
\begin{equation}\label{e10.1}
A=\bigl[A_1\,|\dots|\, A_{k-1}\
\boxed{\ A_k\ }\
A_{k+1}\,|\dots|\,A_t\bigr]
\end{equation}
over $\mathbb F$, partitioned
into vertical strips, among
which one strip $A_k$ is square
and boxed. The number $n_k$ of
rows of $A$ and the number $n_i$
of columns of each strip $A_i$
are nonnegative integers. Let
\begin{equation}\label{e10.8}
B=\bigl[B_1\,|\dots|\, B_{k-1}\
\boxed{\ B_k\ }\
B_{k+1}\,|\dots|\,B_t\bigr]
\end{equation}
be another bangle with the same
sizes of strips and the same $k$
and $t$. We say that the bangles
$A$ and $B$ are
*\!\emph{congruent} or,
respectively, \emph{similar} and
write
\begin{equation}\label{eqeqv}
A\overset{*}{\sim}B
\quad\text{or}\quad
A\overset{s}{\sim}B
\end{equation}
if there exists a nonsingular
upper block-triangular matrix
\begin{equation*}\label{e10.2nn}
S=\begin{bmatrix}
  S_{11}&\dots & S_{1t} \\
  &\ddots&\vdots\\
 0&& S_{tt}
\end{bmatrix}\qquad
(S_{ii}\text{ is } n_i\times
n_i)
\end{equation*}
over $\mathbb F$ such that
\begin{equation*}\label{e10.2}
B=S_{kk}^*AS \quad\text{or}\quad
B=S_{kk}^{-1}AS.
\end{equation*}
\end{definition}

Then
\[
B_k=S_{kk}^*A_kS_{kk}
\quad\text{or}\quad
B_k=S_{kk}^{-1}A_kS_{kk},
\]
this means that the boxed strips
of *congruent/similar bangles
are *congruent/similar. The
following lemma is obvious.

\begin{lemma}\label{s-res1}
Two bangles are
{\rm*}\!congruent/similar if and
only if one reduces to the other
by a sequence of the following
transformations:
\begin{itemize}
  \item[\rm(a)]
Any transformation with rows of
the whole matrix, and then the
{\rm*}\!congruent/similar
transformation with columns of
the boxed strip $($this
transformation reduces
\eqref{e10.1} to
\[
\bigl[\,EA_1\,|\dots| EA_{k-1}\
\boxed{\ EA_kE^*\ }\
EA_{k+1}\,|\dots|\,EA_t\bigr]
\]
or, respectively,
\[
\bigl[\,EA_1\,|\dots| EA_{k-1}\
\boxed{\ EA_kE^{-1}\ }\
EA_{k+1}\,|\dots|\,EA_t\bigr]
\]
with a nonsingular $E)$.
  \item[\rm(b)]
Any transformation with columns
of an unboxed strip.
  \item[\rm(c)]
Addition of a linear combination
of columns of the $i^{\rm th}$
strip to a column of the $j^{\rm
th}$ strip if
$i<j$.\hfill$\square$
\end{itemize}
\end{lemma}

Note that the canonical form
problem for matrices of forms
\eqref{0.0} and mappings
\eqref{m.4} is the canonical
form problem for bangles
\eqref{e10.1} with two strips.
But applying our algorithm to
bangles with two strips we can
produce bangles with three
strips (see Section
\ref{sub-r2}); so we consider
bangles with an arbitrary number
of strips.
\medskip

The paper is organized as
follows. In Section \ref{s-a} we
formulate our main theorem about
the existence of a regularizing
decomposition of a bangle. In
Sections \ref{s-red} and
\ref{s-redz} we construct
regularizing decompositions of
bangles with respect to
*congruence and similarity. In
Section \ref{s-pr} we use these
decompositions to prove the main
theorem and Theorems \ref{t0.01}
and \ref{t.m}.

\section{Bangles}
\label{s-a}

In this section, we formulate
our main theorem, which reduces
the canonical form problem for
bangles up to
*congruence/similarity to the
canonical form problem for
nonsingular matrices up to
*congruence/similarity, and
solves it for complex bangles.

By the \emph{block-direct sum}
of two bangles \eqref{e10.1} and
\eqref{e10.8} with the same
number of strips and the same
position of the boxed strip, we
mean the bangle
\[
A\uplus B:=
\left[\!\begin{array}{c|ccc|c}
\begin{matrix}
  A_1&0 \\
  0&B_1
\end{matrix}&
\begin{matrix}
 \dots\!\!\! \\ \dots\!\!\!
\end{matrix}&
\boxed{\;
\begin{matrix}
    A_k&0 \\
  0&B_k
\end{matrix}\;}&
\begin{matrix}
 \!\!\dots \\ \!\!\dots
\end{matrix}&
\begin{matrix}
  A_t&0 \\
  0&B_t
\end{matrix}
\end{array}\!\right].
\]

\begin{definition}\label{defin}
A \emph{regularizing
decomposition} of a bangle
\begin{equation*}%\label{e10.1q}
A=\bigl[A_1\,|\dots|\, A_{k-1}\
\boxed{\ A_k\ }\
A_{k+1}\,|\dots|\,A_t\bigr]
\end{equation*}
over a field or skew field\/
$\mathbb F$ with respect to
{\rm*\!}congruence/similarity is
a bangle $\Sigma_A$ satisfying
two conditions:
\begin{itemize}
  \item[\rm(i)]
$\Sigma_A$ is
{\rm*\!}congruent/similar to
$A$, and
  \item[\rm(ii)]
$\Sigma_A$ is the block-direct
sum of
\begin{itemize}
  \item
its \emph{regular part}
\begin{equation}\label{e10.3}
\bigl[0_{p0}\,|\dots| \,0_{p0}\
\boxed{\;K\;}\ \,
0_{p0}\,|\dots|\,0_{p0}\bigr],
\qquad \text{$K$ is
nonsingular,}
\end{equation}
  \item
and its \emph{singular part}
being a block-direct sum of
matrices of the form
\begin{gather}\label{e10.5a}
\bigl[0_{q0}\,|\dots| \,0_{q0}\
\boxed{\;J_q(0)\;}\ \,
0_{q0}\,|\dots|\,0_{q0}\bigr],
       \\ \label{e10.5c}
\bigl[\dots|\,E_q\,|\dots \,
\boxed{\;J_q(0)\;}\,
\dots\bigr],
    \qquad
\bigl[\dots\,
\boxed{\;J_q(0)\;}\,
\dots|\,E_q\,|\dots\bigr],
\end{gather}
in which $E_q$ is defined in
\eqref{bip1s} and the dots
denote sequences of strips
$0_{q0}$.
\end{itemize}
Both the regular and the
singular parts may have size
$0$-by-$0$.
\end{itemize}
\end{definition}

The following theorem
generalizes Theorems \ref{t0.01}
and \ref{t.m}.

\begin{theorem}\label{l10.1}
{\rm(a)} Over a field or skew
field $\mathbb F$, any bangle
$A$ possesses regularizing
decompositions for
{\rm*\!}congruence and for
similarity, their regular parts
are determined by $A$ uniquely
up to {\rm*\!}congruence and,
respectively, similarity, and
their singular parts are
determined by $A$ uniquely up to
permutation of summands.

{\rm(b)} If\/ $\mathbb F=\mathbb
C$ and $\Sigma_A$ is a
regularizing decomposition of a
bangle $A$ for to
{\rm*\!}congruence, then its
regular part \eqref{e10.3} is
{\rm *\!}congruent to the
block-direct sum
\[
\biguplus_i\,
\bigl[0_{p_i0}\,|\dots|
\,0_{p_i0}\ \boxed{\;K_i\;}\ \,
0_{p_i0}\,|\dots|\,0_{p_i0}\bigr],
\]
in which $K_1\oplus\dots\oplus
K_s$ is the canonical form of
$K$ defined in Theorem
\ref{bilin} and each $K_i$ is
$p_i$-by-$p_i$. Replacing in
$\Sigma_A$ the regular part by
this block-direct sum, we obtain
a canonical form of $A$ for {\rm
*\!}congruence $($in particular,
for congruence$)$ since the
obtained bangle is {\rm
*\!}congruent to $A$ and is
determined by $A$ uniquely up to
permutation of summands.

{\rm(c)} If\/ $\mathbb F=\mathbb
C$ and $\Sigma_A$ is a
regularizing decomposition of a
bangle $A$ for similarity, then
its regular part is similar to a
block-direct sum of matrices of
the form
\[
\bigl[0_{q0}\,|\dots| \,0_{q0}\
\boxed{\;J_q(\lambda)\;}\ \,
0_{q0}\,|\dots|\,0_{q0}\bigr],
\qquad\lambda\ne 0.
\]
Replacing in $\Sigma_A$ the
regular part by this
block-direct sum, we obtain a
canonical form of $A$ for
similarity since the obtained
bangle is similar to $A$ and is
determined by $A$ uniquely up to
permutation of summands.
\end{theorem}

Note that for bangles with
respect to similarity this
theorem can be deduced from the
canonical form problem solved in
\cite[\S\,2]{naz}.

\section{Regularization
for *congruence}\label{s-red}

We give an algorithm that for
every bangle over a field or
skew field $\mathbb F$
constructs its regularizing
decomposition for *congruence.
If $\mathbb F=\mathbb C$, then
we can improve the numerical
stability of this algorithm
using only unitary
transformations. The algorithm
is the alternating sequence of
left-hand and right-hand
reductions, which we define in
Sections \ref{sub-r1} and
\ref{sub-r2}.

\subsection{Left-hand reduction for
*congruence}\label{sub-r1}

Let
\begin{equation}\label{erre}
A=\bigl[\,A_1\,|\dots|\,
A_{k-1}\ \boxed{\ A_k\ }\
A_{k+1}\,|\dots|\,A_t\,\bigr]
\end{equation}
be a bangle over $\mathbb F$.
Producing *congruence
transformations (a)--(c) from
Lemma \ref{s-res1} with $A$, we
can reduce its submatrix
$[\,A_1| A_2\,|
\dots|\,A_{k-1}\,]$ by the
following transformations:
\begin{itemize}
  \item[(a$'$)]
arbitrary transformations of
rows;
  \item[(b$'$)]
arbitrary transformations of
columns within any vertical
strip $A_i$;
  \item[(c$'$)]
addition of a linear combination
of columns of the $i^{\rm th}$
strip to a column of the $j^{\rm
th}$ strip if $i<j$.
\end{itemize}

First we reduce $[\,A_1| A_2\,|
\dots|\,A_{k-1}\,]$ to the form
\begin{equation}\label{eqqe}
 \left[\begin{array}{cc|c|c|c}
 0&0&A_{12}&\dots&A_{1,k-1}\\
 0&I&A_{22}&\dots&A_{2,k-1}
\end{array}\right]
\end{equation}
using transformations (b$'$)
with $A_1$ and (a$'$), then make
zero $A_{22},\dots,A_{2,k-1}$ by
transformations (c$'$).
Transforming analogously the
submatrix $[\,A_{12}\,|
\dots|\,A_{1,k-1}\,]$, we reduce
\eqref{eqqe} to the form
\[
 \left[\begin{array}{cc|cc|c|c|c}
 0&0 & 0&0 & B_{3}&\dots&B_{k-1}\\
 0&0 & 0&I & 0&\dots&0\\
 0&I & 0&0 & 0&\dots&0
\end{array}\right];
\]
and so on. Repeat this process
until obtain
\begin{equation}\label{lab}
\left[\begin{array}{c|c|c|c}
\begin{matrix}
0&0\\0&0\\\vdots&\vdots
\\0&0\\0&\!I_{r_k}\!
\end{matrix}&
\begin{matrix}
0&0\\0&0\\\vdots&\vdots\\
0&\!I_{r_{k-1}}\!\!\\0&0
\end{matrix}
&
\begin{matrix}
 \dots \\ \dots\\ \phantom{\vdots}\\
  \dots \\ \dots
\end{matrix}&
\begin{matrix}
0&0\\0&\!I_{r_2}\!\!\!\\\vdots&\vdots
\\0&0\\0&0
\end{matrix}
\end{array}\right],\qquad r_2\ge
0,\ \dots,\ r_{k}\ge 0,
\end{equation}
and extend the obtained
partition into horizontal strips
to the whole bangle
\eqref{erre}. Make zero all
horizontal strips of the blocks
$A_k,\dots,A_t$ except for the
first strip and obtain
\begin{equation}\label{el.2}
{\cal L}_k(M):=
\left[\!\begin{array}{c|c|ccc|c|c}
\begin{matrix}
0&0\\0&0\\\vdots&\vdots\\
0&\!I_{r_{k}}\!\!
\end{matrix}
&
\begin{matrix}
 \dots \\ \dots\\ \phantom{\vdots}\\
  \dots
\end{matrix}&
\begin{matrix}
0&0\\0&\!I_{r_{2}}\!\!\!\\\vdots&
\vdots\\0&0
\end{matrix}
 &
\!\boxed{ \
\begin{matrix}
   M_1&M_2&\dots &M_k\\
0&0_{r_{2}}&\dots &0\\
\vdots&\vdots&\ddots& \vdots\\
0&0&\dots &0_{r_{k}}
\end{matrix}\;}\!&\!\!
\begin{matrix}
M_{k+1}\\0\\\vdots\\0
\end{matrix}
&
\begin{matrix}
 \dots \\ \dots\\ \phantom{\vdots}\\
  \dots
\end{matrix}&
\begin{matrix}
M_t\\0\\\vdots\\0
\end{matrix}
\end{array}\!\right]
\end{equation}
(we have divided the boxed block
$A_k$ into $k$ vertical strips
conformally to its partition
into horizontal strips) for some
\begin{equation}\label{el.3}
M =\bigl[\,\boxed{\;M_1\;}\
M_2\,|\dots|\,M_{t}\,\bigr]=:L(A).
\end{equation}
Clearly, $r_2,\dots,r_{k}$ are
uniquely determined by $A$.

\begin{definition}
We say that a bangle $A$ reduces
to a bangle $B$ by
\emph{admissible permutations}
and write
\begin{equation*}\label{el.4}
A\overset{p}{\sim}B
\end{equation*}
if $A$ reduces to $B$ by a
sequence of the following
transformations:
\begin{itemize}
  \item
permutation of rows of the whole
matrix and then the same
permutation of columns of the
boxed strip,
  \item
permutation of columns in an
unboxed strip.
\end{itemize}
\end{definition}

Clearly,
\[
A\overset{p}{\sim}B
 \qquad\Longrightarrow\qquad
A\overset{s}{\sim}B\ \text{ and
}\ A\overset{*}{\sim}B
\]
(in the notation \eqref{eqeqv}).

\begin{lemma}\label{lem_r1}
{\rm(a)} The equivalence
\begin{equation}\label{el.5}
{\cal L}_k(M)\overset{*}{\sim}
{\cal L}_k(N)
 \quad\Longleftrightarrow\quad
 M\overset{*}{\sim} N
\end{equation}
holds for all
\begin{equation*}
M=\bigl[\,\boxed{\;M_1\;}\
M_2\,|\dots|\,M_{t}\,\bigr],
\qquad
N=\bigl[\,\boxed{\;N_1\;}\
N_2\,|\dots|\,N_{t}\,\bigr],
\end{equation*}
and each $k\le t$.

{\rm(b)} If\/ $\mathbb F=
\mathbb C$, then for every
bangle $A$ we can find
\eqref{el.2} using only unitary
transformations.
\end{lemma}

\begin{proof}
(a) The equivalence \eqref{el.5}
is trivial if $k=1$. Let $k\ge
2$. Reasoning by induction on
$k$, we assume that
\begin{equation}\label{el.5a}
{\cal
L}_{k-1}(M)\overset{*}{\sim}
{\cal L}_{k-1}(N)
 \quad\Longleftrightarrow\quad
 M\overset{*}{\sim} N
\end{equation}
and prove the equivalence
\eqref{el.5} as follows.

\begin{itemize}
  \item[$(\Rightarrow)$]
Suppose ${\cal
L}_k(M)\overset{*}{\sim}{\cal
L}_k(N)$, that is,
\begin{equation}\label{el.6a}
S_{kk}^*{\cal L}_k(M)S={\cal
L}_k(N)
\end{equation}
for some nonsingular
\begin{equation}\label{mdid}
S=\begin{bmatrix}
  S_{11}&\dots & S_{1t} \\
  &\ddots&\vdots\\
  0&& S_{tt}
\end{bmatrix}.
\end{equation}
Since both ${\cal L}_k(M)$ and
${\cal L}_k(N)$ have the same
first vertical strip
\[
\begin{bmatrix}
  0&0 \\
  0&I_{r_k}
\end{bmatrix}
\]
(we join its zero horizontal
strips), by \eqref{el.6a} we
have
\[
S_{kk}^*\begin{bmatrix}
  0&0 \\
  0&I_{r_k}
\end{bmatrix}
S_{11}=\begin{bmatrix}
  0&0 \\
  0&I_{r_k}
\end{bmatrix}
\]
and so $S_{kk}$ has the form
\begin{equation}\label{el.6d}
S_{kk}=\begin{bmatrix}
  P_1 & P_2 \\
  0 & P_3
\end{bmatrix}.
\end{equation}
Let
\[
R:=\begin{bmatrix}
  S_{22}&\dots & S_{2t} \\
  &\ddots&\vdots\\
 0 && S_{tt}
\end{bmatrix}
\]
be a submatrix of \eqref{mdid}
with $S_{kk}$ of the form
\eqref{el.6d}. Due to
\eqref{el.6a},
\begin{equation}\label{el.6aq}
P_1^*{\cal L}_{k-1}(M)R={\cal
L}_{k-1}(N).
\end{equation}
 So ${\cal
L}_{k-1}(M)\overset{*}{\sim}
{\cal L}_{k-1}(N)$, and by
\eqref{el.5a}
$M\overset{*}{\sim} N$.

  \item[$(\Longleftarrow)$]
Suppose $M\overset{*}{\sim}N$.
By \eqref{el.5a}, ${\cal
L}_{k-1}(M)\overset{*}{\sim}
{\cal L}_{k-1}(N)$, this ensures
\begin{equation*}\label{el.6ad}
P_{kk}^*{\cal L}_{k-1}(M)P={\cal
L}_{k-1}(N)
\end{equation*}
for some nonsingular

\begin{equation*}\label{mdi}
P=\begin{bmatrix}
  P_{11}&\dots & P_{1t} \\
  &\ddots&\vdots\\
  0&& P_{tt}
\end{bmatrix}.
\end{equation*}
Denote by $B_i$ and $C_i$ the
strips of ${\cal L}_{k-1}(M)$
and ${\cal L}_{k-1}(N)$:
\[
{\cal
L}_{k-1}(M)=\bigl[\,B_1\,|\dots\,|\,
B_{k-1}\ \boxed{\ B_k\ }\
B_{k+1}\,|\,\dots|\,B_t\,\bigr],\]
\[{\cal
L}_{k-1}(N)=\bigl[\,C_1\,|\dots\,|\,
C_{k-1}\ \boxed{\ C_k\ }\
C_{k+1}\,|\,\dots|\,C_t\,\bigr].
\]
Then
\[
{\cal L}_k(M)=
\left[\!\begin{array}{c|c|c|ccc|c|c}
\begin{matrix}
  0&0 \\
  0&I_{r_k}
\end{matrix}&\begin{matrix}
  B_1 \\
  0
\end{matrix}&
\begin{matrix}
 \dots \\ \dots
\end{matrix}&\begin{matrix}
  B_{k-1}\!\!\\0\!\!
\end{matrix}
& \boxed{\;
\begin{matrix}
    B_k&B_{k+1} \\
  0&0
\end{matrix}\;}
&\begin{matrix}
 \!\! B_{k+2}\\\!\! 0
\end{matrix}
&
\begin{matrix}
 \dots \\ \dots
\end{matrix}&
\begin{matrix}
  B_t\\0
\end{matrix}
\end{array}\!\right]
\]
and by \eqref{el.6aq}
\begin{multline*}\label{e10.11}
{\cal L}_k(M)\overset{*}{\sim}
\begin{bmatrix}
  P_{kk} & P_{k,k+1} \\
  0 & P_{k+1,k+1}
\end{bmatrix}^*{\cal L}_k(M)
\begin{bmatrix}
\,
\begin{bmatrix}
  I&0\\
  0&(P_{k+1,k+1}^*)^{-1}
\end{bmatrix}
  &0\\
0 & P
\end{bmatrix}
  \\ =
\left[\begin{array}{cc|c|c|c}
 0&0&C_1& \dots&C_t\\
 0&I_{r_k}&C'_1& \dots&C'_t
\end{array}\right]\overset{*}{\sim}
{\cal L}_k(N),
\end{multline*}
where $C'_1,\dots,C'_{t}$ are
some matrices.
\end{itemize}

This proves \eqref{el.5}. Let us
give an alternative proof of
\eqref{el.5} using *congruence
transformations (a)--(c) from
Lemma \ref{s-res1}. Due to that
lemma, it suffices to show that
those transformations (a)--(c)
with \eqref{el.2} that preserve
all of its blocks except for
$M_1,\dots,M_{t}$ produce all
transformations (a)--(c) with
\eqref{el.3}.

\begin{itemize}
  \item
We can add a column of $M_i$ to
a column of $M_j$ if $i<j$.
Indeed, in the case $j\le k$
this is a column-transformation
within the boxed block of ${\cal
L}_k(M)$, and so we must produce
the *congruent
row-transformation---add the
corresponding row of the $i^{\rm
th}$ horizontal strip of
\eqref{el.2} to the row of the
$j^{\rm th}$ horizontal strip.
This spoils zero blocks of the
$j^{\rm th}$ horizontal strip,
but they are repaired by
additions of columns of
$I_{r_{j}}$.

  \item
We can also make arbitrary
elementary transformations with
columns of $M_i$ if $i\ne 1$: in
the case $i\le k$ these
transformations spoil $I_{r_i}$
but it is restored by
transformations with its
columns.

\end{itemize}
\medskip

(b) Let $\mathbb F=\mathbb C$.
We must prove that if
\[
A=\bigl[A_1\,|\dots|\, A_{k-1}\
\boxed{\ A_k\ }\
A_{k+1}\,|\dots|\,A_t\bigr]
\]
is reduced to \eqref{el.3} by
the algorithm from this section,
then $r_2,\dots,r_k$ and
$M_1,\dots,M_t$ can be found
using only unitary
transformations with $A$. By
unitary column-transformations
within vertical strips
$A_1,\dots,A_{k-1}$ of $A$ and
by unitary
row-trans\-formations, we
sequentially reduce its
submatrix $[\,A_1| A_2\,|
\dots|\,A_{k-1}\,]$  to the form
\begin{equation*}\label{el.1aa}
\left[\begin{array}{c|c|c|c}
\begin{matrix}
0&0\\0&0\\0&0\\\vdots&\vdots
\\0&0\\0&\!H_{r_k}\!
\end{matrix}&
\begin{matrix}
0&0\\0&0\\0&0\\\vdots&\vdots\\
0&\!H_{r_{k-1}}\!\!\\ *&*
\end{matrix}
&
\begin{matrix}
 \dots \\\dots \\ \dots\\ \phantom{\vdots}\\
  \dots \\ \dots
\end{matrix}&
\begin{matrix}
0&0\\0&\!H_{r_2}\!\!\!\\
*&*\\\vdots&\vdots
\\ *&*\\ *&*
\end{matrix}
\end{array}\right],
\end{equation*}
where each $H_{r_i}$ is a
nonsingular $r_i$-by-$r_i$ block
and all $*$'s are unspecified
blocks (this reduction was
studied thoroughly in
\cite[Section 4]{ser}). The
matrix $A$ takes the form
\begin{equation}\label{lwl}
\left[\begin{array}{c|c|ccc|c|c}
\begin{matrix}
0&0\\0&0\\\vdots&\vdots\\
0&\!H_{r_{k}}\!\!
\end{matrix}
&
\begin{matrix}
 \dots \\ \dots\\ \phantom{\vdots}\\
  \dots
\end{matrix}&
\begin{matrix}
0&0\\0&\!H_{r_{2}}\!\!\!\\\vdots&
\vdots\\ *&*
\end{matrix}
 &
\boxed{ \
\begin{matrix}
   M_1&M_2&\dots &M_k\\
*&*_{r_{2}}&\dots &*\\
\vdots&\vdots&\ddots& \vdots\\
*&*&\dots &*_{r_{k}}
\end{matrix}\;}&\!\!
\begin{matrix}
M_{k+1}\\ *\\\vdots\\ *
\end{matrix}
&
\begin{matrix}
 \dots \\ \dots\\ \phantom{\vdots}\\
  \dots
\end{matrix}&
\begin{matrix}
M_t\\ *\\\vdots\\ *
\end{matrix}
\end{array}\right],
\end{equation}
in which
$*_{r_{2}},\dots,*_{r_{k}}$ are
$r_2\times r_2,\dots, r_k\times
r_k$ matrices. Replacing
$H_{r_{2}},\dots,H_{r_{k}}$ by
the identity matrices of the
same sizes and all $*$'s by the
zero matrices, we obtain
\eqref{el.3} because \eqref{lwl}
can be reduced to \eqref{el.3}
by those transformations
(a)--(c) from Lemma \ref{s-res1}
that preserve $r_2,\dots,r_k$
and $M_1,\dots,M_t$.
\end{proof}

\subsection{Right-hand reduction for
*congruence}\label{sub-r2}

Let
\begin{equation}\label{el.31}
A=\bigl[\:\boxed{\ A_1\ }\
A_{2}\,|\dots|\,A_t\,\bigr]
\end{equation}
be a bangle over a field or skew
field $\mathbb F$.

First we reduce $A$ by
*congruence transformations
\begin{equation}\label{nnn}
\bigl[\:\boxed{\ SA_1S^*\ }\
SA_{2}\,|\dots|\,SA_t\,\bigr],\qquad
\text{$S$ is nonsingular,}
\end{equation}
to the form
\begin{equation}\label{el.32}
\left[\begin{array}{cc|c|c}
\boxed{\;
\begin{matrix}
   0_d&0\\
  B_1'&B_{2}'
\end{matrix}\;}\!&
\begin{matrix}
  B_{3}\\B_3'
\end{matrix}&
\begin{matrix}
 \dots \\ \dots
\end{matrix}&
\begin{matrix}
  B_{t+1}\\B_{t+1}'
\end{matrix}
\end{array}\right],
\end{equation}
in which the rows of $[B_1'\
B_{2}']$ are linearly
independent and $B_{2}'$ is
square.

Then we make zero
$B_{3}',\dots,B_{t+1}'$ adding
columns of $B_1'$ and $B_2'$,
and as in \eqref{lab}
sequentially reduce $[\,B_3|
B_4\,| \dots|\,B_{t+1}\,]$ to
the form
\[
\left[\begin{array}{c|c|c|c}
\begin{matrix}
0&0\\0&0\\\vdots&\vdots
\\0&0\\ 0&\!I_{r_t}\!
\end{matrix}&
\begin{matrix}
0&0\\0&0\\\vdots&\vdots\\
0&I_{r_{t-1}}\\0&0
\end{matrix}
&
\begin{matrix}
 \dots \\ \dots\\ \phantom{\vdots}\\
  \dots \\ \dots
\end{matrix}&
\begin{matrix}
0&0\\0&I_{r_2}\\\vdots&\vdots
\\0&0\\0&0
\end{matrix}
\end{array}\right],
\]
obtaining a partition of the
first horizontal strip of
\eqref{el.32} into $t$
substrips. Conformally divide
the first vertical strip of the
boxed block into $t$ substrips
and obtain
\begin{multline}\label{el.33}
{\cal R}(M)=
 \\
\left[\begin{array}{cc|c|c|c}
\boxed{ \
\begin{matrix}
   0_{r_1}&0&\dots &0&0&0\\
0&0_{r_2}&\dots &0&0&0\\
\vdots&\vdots&\ddots&
\vdots&\vdots&\vdots\\ 0&0&\dots
&0_{r_{t-1}}&0&0\\
 0&0&\dots &0&0_{r_{t}}&0\\
  M_1&M_2&\dots&M_{t-1}&M_{t}&M_{t+1}
\end{matrix}\;}&\!\!
\begin{matrix}
0&0\\0&0\\\vdots&\vdots\\0&0\\0&I_{r_t}\\0&0
\end{matrix}&
\begin{matrix}
0&0\\0&0\\\vdots&\vdots\\
0&I_{r_{t-1}}\\0&0\\0&0
\end{matrix}
&
\begin{matrix}
 \dots \\ \dots\\ \phantom{\vdots}\\
  \dots \\ \dots\\ \dots
\end{matrix}&
\begin{matrix}
0&0\\0&I_{r_2}\\\vdots&\vdots\\0&0\\0&0\\0&0
\end{matrix}
\end{array}\right]
\end{multline}
for some
\begin{equation}\label{el.34}
M =\bigl[\,M_1|\dots|M_{t}\
\boxed{\;M_{t+1}\;}\
\bigr]=:R(A)
\end{equation}
with $M_{t+1}=B_2$.

\begin{lemma}\label{lem_r2}
{\rm(a)} The equivalence
\begin{equation}\label{el.35}
{\cal R}(M)\overset{*}{\sim}
{\cal R}(N)
 \quad\Longleftrightarrow\quad
 M\overset{*}{\sim} N
\end{equation}
holds for all
\begin{equation*}
M=\bigl[\,M_1|\dots|M_{t}\
\boxed{\;M_{t+1}\;}\ \bigr],
\qquad
N=\bigl[\,N_1|\dots|N_{t}\
\boxed{\;N_{t+1}\;}\ \bigr].
\end{equation*}

{\rm(b)} If\/ $\mathbb F=
\mathbb C$, then for every
bangle $A$ of the form
\eqref{el.31} we can find
\eqref{el.33} using only unitary
transformations.
\end{lemma}

\begin{proof}
(a) Let us prove the equivalence
\eqref{el.35} using *congruence
transformations (a)--(c) from
Lemma \ref{s-res1}
(alternatively, one could use
induction on $t$ as in the proof
of Lemma \ref{lem_r1}(a)). Due
to Lemma \ref{s-res1}, it
suffices to show that those
transformations (a)--(c) with
\eqref{el.33} that preserve all
of its blocks except for
$M_1,\dots,M_{t+1}$ produce all
transformations (a)--(c) with
\eqref{el.34}.
\begin{itemize}
  \item
We can add a column of $M_i$ to
a column of $M_j$ if $i<j$; by
the definition of *congruence
transformations we must add the
corresponding row of the $i^{\rm
th}$ horizontal strip of
\eqref{el.33} to the row of the
$j^{\rm th}$ horizontal strip;
although this spoils a zero
block of the $j^{\rm th}$
horizontal strip if $i\ne 1$,
but it can be repaired by
additions of columns of
$I_{r_j}$.
  \item
We can also make arbitrary
elementary transformations with
columns of $M_i$ if $i\le t$:
these transformations spoil
$I_{r_i}$ if $i\ne 1$, but it is
restored by transformations with
its columns.
\end{itemize}

(b) Let $\mathbb F=\mathbb C$.
First we reduce the bangle
\eqref{el.31} by transformations
\eqref{nnn} with unitary $S$ to
the form \eqref{el.32}, in which
the rows of $[B_1'\ B_{2}']$ are
linearly independent and
$B_{2}'$ is square.

Then we sequentially reduce
$[\,B_3| B_4\,|
\dots|\,B_{t+1}\,]$
 by
unitary column-transformations
within vertical strips and by
unitary row-transformations to
the form
\[
\left[\begin{array}{c|c|c|c}
\begin{matrix}
0&0\\0&0\\\vdots&\vdots
\\0&0\\ 0&\!H_{r_t}\!
\end{matrix}&
\begin{matrix}
0&0\\0&0\\\vdots&\vdots\\
0&H_{r_{t-1}}\\ *&*
\end{matrix}
&
\begin{matrix}
 \dots \\ \dots\\ \phantom{\vdots}\\
  \dots \\ \dots
\end{matrix}&
\begin{matrix}
0&0\\0&H_{r_2}\\\vdots&\vdots
\\ *&*\\ *&*
\end{matrix}
\end{array}\right],
\]
where each $H_{r_i}$ is a
nonsingular $r_i$-by-$r_i$ block
and the $*$'s are unspecified
blocks. The matrix $A$ takes the
form
\begin{equation}\label{mkm}
\left[\begin{array}{cc|c|c}
\boxed{ \
\begin{matrix}
   0_{r_1}&0&\dots &0&0\\
0&0_{r_2}&\dots &0&0\\
\vdots&\vdots&\ddots&
\vdots&\vdots\\
 0&0&\dots &0_{r_{t}}&0\\
  M_1&M_2&\dots&M_{t}&M_{t+1}
\end{matrix}\;}&\!\!
\begin{matrix}
0&0\\0&0\\\vdots&\vdots\\0&H_{r_t}\\
* &*
\end{matrix}&
\begin{matrix}
 \dots \\ \dots\\ \phantom{\vdots}\\
  \dots\\ \dots
\end{matrix}&
\begin{matrix}
0&0\\0&H_{r_2}\\\vdots&\vdots\\
*&*\\ *&*
\end{matrix}
\end{array}\right],
\end{equation}
where $M_{t+1}=B_2$. Replacing
$H_{r_{2}},\dots,H_{r_{k}}$ by
the identity matrices of the
same sizes and all $*$'s by the
zero matrices, we obtain
\eqref{el.33} because
\eqref{mkm} can be reduced to
\eqref{el.33} by those
transformations (a)--(c) from
Lemma \ref{s-res1} that preserve
$r_1,\dots,r_t$ and
$M_1,\dots,M_{t+1}$.
\end{proof}

\subsection{Regularization algorithm
for *congruence}\label{sub-r3}

For any bangle
\begin{equation}\label{ooo}
A=\bigl[\,A_1\,|\dots|\,
A_{k-1}\ \boxed{\ A_k\ }\
A_{k+1}\,|\dots|\,A_t\,\bigr]
\end{equation}
over $\mathbb F$, its
regularizing decomposition for
*congruence can be constructed
as follows.

Alternating the left-hand and
the right-hand reductions for
*congruence, we construct the
sequence of bangles
\begin{equation*}
A':=L(A),\ \ A'':=R(A'),\ \
A''':=L(A''),\ \
A'''':=R(A'''),\dots
\end{equation*}
until obtain
\begin{equation}\label{el.63}
A^{(n)}=\bigl[\ \boxed{\ K\ }\
0_{p0}\,|\dots|\,0_{p0}\,\bigr]
\quad\text{or}\quad
A^{(n)}=\bigl[\,0_{p0}\,|\dots|\,
0_{p0}\ \boxed{\ K\ }\ \bigr]
\end{equation}
with a nonsingular $K$.

Producing this reduction, we in
each step have deleted the
reduced parts of $A$; say, in
step 1 we reduced $A$ to the
form \eqref{el.2} and took only
its unreduced part $A'=L(A)$.
Let us repeat the reduction of
\eqref{ooo} preserving all the
reduced parts of $A$:
\begin{itemize}
  \item
In step 1 we transform $A$ to
${\cal L}_k(A')$ of the form
\eqref{el.2}.
  \item
In step 2 we reduce its
subbangle $A'$ to ${\cal
R}(A'')$ preserving the other
blocks of ${\cal L}_k(A')$, and
so on.
\end{itemize}
After $n$ steps, instead of
\eqref{el.63} we obtain some
bangle $\hat{A}$, which is
*congruent to $A$. Due to the
next theorem, $\hat A$ is a
regularizing decomposition of
$A$ up to admissible
permutations of rows and
columns.

\begin{theorem}\label{tel.1}
If $A$ is a bangle over a field
or skew field $\mathbb F$, then
$\hat A$ reduces by admissible
permutations of rows and columns
to a regularizing decomposition
of $A$ for {\rm *}\!congruence.
\end{theorem}

\begin{proof} We give a
constructive proof of this
theorem.

By admissible permutations of
rows and columns, $\hat A$
reduces to a block-direct sum of
the bangle \eqref{e10.3} in
which $K$ is the same as in
\eqref{el.63}, and a bangle $D$
in which each row and each
column contains at most one $1$
and its other entries are zero.
We obtain a regularizing
decomposition of $A$ for
*congruence replacing $D$ in
this block-direct sum by
$\Sigma_D$ from the follows
statement.

\begin{equation}\label{43}
\parbox{25em}
{Let $D$ be a bangle in which
each row and each column
contains at most one $1$ and the
other entries are zero. Then $D$
reduces by admissible
permutations of rows and columns
to a block-direct sum $\Sigma_D$
of bangles of the form
\eqref{e10.5a} and
\eqref{e10.5c}.}
\end{equation}

Let us prove \eqref{43}. By
admissible permutations of rows
and columns of $D$, we reduce
its boxed strip $D_k$ to a
direct sum of singular Jordan
blocks. Then we rearrange
columns in each unboxed strip
such that if its $(i,j)$ and
$(i',j')$ entries are $1$ and
$i<i'$, then $j<j'$. It is easy
to see that the obtained bangle
$\Sigma_D$ is a block-direct sum
of bangles of the form
\eqref{e10.5a} and
\eqref{e10.5c}: each singular
Jordan block $J_p(0)$ in the
decomposition of $D_k$ gives the
summand \eqref{e10.5a} if those
row of $D$ that contains the
last (zero) row of $J_p(0)$ is
zero, and the summand
\eqref{e10.5c} otherwise. The
summands \eqref{e10.5c} with
$p=0$ give zero columns in
unboxed strips of $D$.
\end{proof}

%%%%%%%%%%%%%%%%%%%%%%%%%%%%%%%%5

\section{Regularization
for similarity}\label{s-redz}

We give an algorithm that for
every bangle over a field or
skew field $\mathbb F$
constructs its regularizing
decomposition for similarity. If
$\mathbb F=\mathbb C$, then we
can improve the numerical
stability of this algorithm
using only unitary
transformations.

\subsection{Left-hand reduction for
similarity}\label{sub-r1z}

Let
\begin{equation*}
A=\bigl[\,A_1\,|\dots|\,
A_{k-1}\ \boxed{\ A_k\ }\
A_{k+1}\,|\dots|\,A_t\,\bigr]
\end{equation*}
be a bangle over $\mathbb F$.
Using similarity transformations
with $A$, we can reduce its
submatrix $[\,A_1| A_2\,|
\dots|\,A_{k-1}\,]$ by
transformations (a$'$)--(c$'$)
from Section \ref{sub-r1}. We
reduce this submatrix to the
form
\begin{equation*}
\left[\begin{array}{c|c|c|c}
\begin{matrix}
0&\!I_{r_1}\!\\0&0\\\vdots&\vdots
\\0&0\\0&0
\end{matrix}&
\begin{matrix}
0&0\\0&\!I_{r_2}\!\\\vdots&\vdots\\
0&0\\0&0
\end{matrix}
&
\begin{matrix}
 \dots \\ \dots\\ \phantom{\vdots}\\
  \dots \\ \dots
\end{matrix}&
\begin{matrix}
0&0\\0&0\\\vdots&\vdots
\\0&\!I_{r_{k-1}}\!\!\!\\0&0
\end{matrix}
\end{array}\right],\qquad r_1\ge
0,\ \dots,\ r_{k-1}\ge 0,
\end{equation*}
and obtain a partition of the
bangle $A$ into $k$ horizontal
strips. Then we divide the boxed
block $A_k$ into $k$ vertical
substrips of the same sizes,
make zero all horizontal strips
in the blocks $A_k,\dots,A_t$
except for the last strip, and
obtain
\begin{multline}\label{el.2z}
{\cal L}_k(M)=
  \\
\left[\begin{array}{c|c|ccc|c|c}
\begin{matrix}
0&\!I_{r_{1}}\!\!\\\vdots&\vdots\\0&0\\
0&0
\end{matrix}
&
\begin{matrix}
 \dots \\ \phantom{\vdots} \\ \dots\\
  \dots
\end{matrix}&
\begin{matrix}
0&0\\\vdots&
\vdots\\0&\!I_{r_{k-1}}\!\!\!\\0&0
\end{matrix}
 &
\boxed{ \
\begin{matrix}
0_{r_{1}}&\dots &0&0
\\
\vdots&\ddots&\vdots& \vdots\\
0&\dots&0_{r_{k-1}} &0\\
  M_1&\dots &M_{k-1}&M_k
\end{matrix}\;}&\!\!
\begin{matrix}
0\\\vdots\\0\\M_{k+1}
\end{matrix}
&
\begin{matrix}
 \dots \\  \phantom{\vdots}\\\dots\\
  \dots
\end{matrix}&
\begin{matrix}
0\\\vdots\\0\\M_t
\end{matrix}
\end{array}\right]
\end{multline}
for some
\begin{equation}\label{el.3z}
M =\bigl[\,M_1\,|\dots|\,
M_{k-1}\ \boxed{\ M_k\ }\
M_{k+1}\,|\dots|\,M_t\,\bigr]=:L(A).
\end{equation}

\begin{lemma}\label{lem_r1z}
{\rm(a)} The equivalence
\begin{equation*}
{\cal L}_k(M)\overset{s}{\sim}
{\cal L}_k(N)
 \quad\Longleftrightarrow\quad
 M\overset{s}{\sim} N
\end{equation*}
holds for all
\begin{align*}
M &=\bigl[\,M_1\,|\dots|\,
M_{k-1}\ \boxed{\ M_k\ }\
M_{k+1}\,|\dots|\,M_t\,\bigr],\\
N & =\bigl[\,N_1\,|\dots|\,
N_{k-1}\ \boxed{\ N_k\ }\
N_{k+1}\,|\dots|\,N_t\,\bigr].
\end{align*}

{\rm(b)} If\/ $\mathbb F=
\mathbb C$, then for every
bangle $A$ we can find
\eqref{el.2z} using only unitary
transformations.
\end{lemma}

\begin{proof}
(a) This statement follows from
Lemma \ref{s-res1} since those
transformations (a)--(c) with
\eqref{el.2z} that preserve all
of its blocks except for
$M_1,\dots,M_{t}$ produce all
transformations (a)--(c) with
\eqref{el.3z}. For example, we
can add a column of $M_i$ to a
column of $M_j$ if $i<j$:
although in the case $j\le k$ we
must subtract the corresponding
row of the $j^{\rm th}$
horizontal strip of
\eqref{el.2z} from the row of
the $i^{\rm th}$ horizontal
strip, and this may spoil zero
blocks of the $i^{\rm th}$
horizontal strip, but they are
repaired by additions of columns
of $I_{r_{i}}$.
\medskip

(b) Let $\mathbb F=\mathbb C$.
By unitary
column-transformations within
vertical strips of $A$ and by
unitary row-trans\-formations,
we sequentially reduce its
submatrix $[\,A_1| A_2\,|
\dots|\,A_{k-1}\,]$ to the form
\[
\left[\begin{array}{c|c|c|c}
\begin{matrix}
0&\!H_{r_1}\!\\0&0\\\vdots&\vdots
\\0&0\\0&0
\end{matrix}&
\begin{matrix}
*&*\\0&\!H_{r_2}\!\\\vdots&\vdots\\
0&0\\0&0
\end{matrix}
&
\begin{matrix}
 \dots \\ \dots\\ \phantom{\vdots}\\
  \dots \\ \dots
\end{matrix}&
\begin{matrix}
*&*\\ *&*\\\vdots&\vdots
\\0&\!H_{r_{k-1}}\!\!\!\\0&0
\end{matrix}
\end{array}\right],
\]
where each $H_{r_i}$ is a
nonsingular $r_i$-by-$r_i$ block
and all $*$'s are unspecified
blocks. The matrix $A$ takes the
form
\begin{equation}\label{el.6qz}
\left[\begin{array}{c|c|ccc|c|c}
\begin{matrix}
0&\!H_{r_{1}}\!\!\\\vdots&\vdots\\0&0\\
0&0
\end{matrix}
&
\begin{matrix}
 \dots \\ \phantom{\vdots} \\ \dots\\
  \dots
\end{matrix}&
\begin{matrix}
*&*\\\vdots&
\vdots\\0&\!H_{r_{k-1}}\!\!\!\\0&0
\end{matrix}
 &
\boxed{ \
\begin{matrix}
*_{r_{1}}&\dots &*&*
\\
\vdots&\ddots&\vdots& \vdots\\
*&\dots&*_{r_{k-1}} &*\\
  M_1&\dots &M_{k-1}&M_k
\end{matrix}\;}&\!\!
\begin{matrix}
*\\\vdots\\ *\\M_{k+1}
\end{matrix}
&
\begin{matrix}
 \dots \\  \phantom{\vdots}\\\dots\\
  \dots
\end{matrix}&
\begin{matrix}
*\\\vdots\\ *\\M_t
\end{matrix}
\end{array}\right],
\end{equation}
in which
$*_{r_1},\dots,*_{r_{k-1}}$ are
$r_1\times r_1,\dots,
r_{k-1}\times r_{k-1}$ matrices.
Replacing
$H_{r_1},\dots,H_{r_{k-1}}$ by
the identity matrices of the
same sizes and all $*$'s by the
zero matrices, we obtain
\eqref{el.2z} since
\eqref{el.6qz} reduces to
\eqref{el.2z} by those
transformations (a)--(c) from
Lemma \ref{s-res1} that preserve
$r_1,\dots,r_{k-1}$,
$M_1,\dots,M_t$.
\end{proof}

\subsection{Right-hand reduction for
similarity}\label{sub-r2z}

Let
\begin{equation}\label{el.31z}
A=\bigl[\:\boxed{\ A_1\ }\
A_{2}\,|\dots|\,A_t\,\bigr]
\end{equation}
be a bangle over $\mathbb F$.

First we reduce $A$ by
similarity transformations
\begin{equation}\label{nn}
\bigl[\:\boxed{\ SA_1S^{-1}\ }\
SA_{2}\,|\dots|\,SA_t\,\bigr],\qquad
\text{$S$ is nonsingular,}
\end{equation}
to the form
\begin{equation}\label{kok}
\left[\begin{array}{cc|c|c}
\boxed{\;
\begin{matrix}
     B_1&B_{2}\\0&0
\end{matrix}\;}\!&
\begin{matrix}
B_3\\ B_{3}'
\end{matrix}&
\begin{matrix}
 \dots \\ \dots
\end{matrix}&
\begin{matrix}
B_{t+1}\\B_{t+1}'\\
\end{matrix}
\end{array}\right],
\end{equation}
in which the rows of $[B_1\
B_{2}]$ are linearly independent
and $B_{1}$ is square.

Then we make zero $B_3,\,
\dots,\,B_{t+1}$ adding columns
of $B_1$ and $B_2$, and
sequentially reduce $[\,B_3'\,|
\dots|\,B_{t+1}'\,]$ to the form
\[
\left[\begin{array}{c|c|c|c}
\begin{matrix}
0&\!I_{r_2}\!\\0&0\\\vdots&\vdots
\\0&0\\ 0&0
\end{matrix}&
\begin{matrix}
0&0\\0&I_{r_{3}}\\\vdots&\vdots\\
0&0\\0&0
\end{matrix}
&
\begin{matrix}
 \dots \\ \dots\\ \phantom{\vdots}\\
  \dots \\ \dots
\end{matrix}&
\begin{matrix}
0&0\\0&0\\\vdots&\vdots
\\0&\!I_{r_t}\!\\0&0
\end{matrix}
\end{array}\right].
\]

The matrix $A$ transforms to
\begin{multline}\label{el.33z}
{\cal R}(M)=
 \\
\left[\begin{array}{cc|c|c|c}
\boxed{ \
\begin{matrix}
  M_1&M_2&M_3&\dots&M_{t}&M_{t+1}\\
0&0_{r_2}&0&\dots &0&0\\
 0&0&0_{r_3}&\dots &0&0\\
\vdots&\vdots&\vdots&
\ddots&\vdots&\vdots\\
 0&0 &0&\dots&0_{r_{t}}&0\\
0&0&0&\dots &0&0_{r_{t+1}}\\
\end{matrix}\;}&\!\!
\begin{matrix}
0&0\\0&I_{r_2}\\0&0\\
\vdots&\vdots\\0&0\\0&0
\end{matrix}&
\begin{matrix}
0&0\\0&0\\0&I_{r_3}\\\vdots&\vdots\\
0&0\\0&0
\end{matrix}
&
\begin{matrix}
\dots \\ \dots \\ \dots\\
\phantom{\vdots}\\
   \dots\\ \dots
\end{matrix}&
\begin{matrix}
0&0\\0&0\\0&0\\\vdots&\vdots
\\0&I_{r_{t}}\\0&0
\end{matrix}
\end{array}\right],
\end{multline}
for some
\begin{equation}\label{el.34z}
M =\bigl[\:\boxed{\ M_1\ }\
M_{2}\,|\dots|\,M_{t+1}\,\bigr]=:R(A)
\end{equation}
with $M_{1}=B_1$.
\begin{lemma}\label{lem_r2z}
{\rm(a)} The equivalence
\begin{equation*}
{\cal R}(M)\overset{s}{\sim}
{\cal R}(N)
 \quad\Longleftrightarrow\quad
 M\overset{s}{\sim} N
\end{equation*}
holds for all
\begin{equation*}\label{el.36z}
M=\bigl[\:\boxed{\ M_1\ }\
M_{2}\,|\dots|\,M_{t+1}\,\bigr],
\qquad N=\bigl[\:\boxed{\ N_1\
}\
N_{2}\,|\dots|\,N_{t+1}\,\bigr].
\end{equation*}

{\rm(b)} If\/ $\mathbb F=
\mathbb C$, then for every
bangle $A$ of the form
\eqref{el.31z} we can find
\eqref{el.33z} using only
unitary transformations.
\end{lemma}

\begin{proof}
(a) It is easy to show that
those transformations (a)--(c)
from Lemma \ref{s-res1} with
\eqref{el.33z} that preserve all
of its blocks except for
$M_1,\dots,M_{t+1}$ produce all
transformations (a)--(c) with
\eqref{el.34z}. Say, we can add
a column of $M_i$ to a column of
$M_j$ if $i<j$: although we must
subtract the corresponding row
of the $j^{\rm th}$ horizontal
strip of \eqref{el.33z} from the
row of the $i^{\rm th}$
horizontal strip, and this
spoils zero blocks of the
$i^{\rm th}$ horizontal strip if
$j\ne t+1$, but they are
repaired by additions of columns
of $I_{r_i}$.

(b) Let $\mathbb F=\mathbb C$.
First we reduce $A$ by
transformations \eqref{nn} with
unitary $S$ to the form
\eqref{kok}, in which the rows
of $[B_1\ B_{2}]$ are linearly
independent and $B_{1}$ is
square.

Then we sequentially reduce its
submatrix $[\,B_3'| B_4'\,|
\dots|\,B_{t+1}'\,]$ by unitary
column-transformations within
vertical strips and by unitary
row-transformations to the form
\[
\left[\begin{array}{c|c|c|c}
\begin{matrix}
0&\!H_{r_2}\!\\0&0\\\vdots&\vdots
\\0&0\\ 0&0
\end{matrix}&
\begin{matrix}
*&*\\0&H_{r_{3}}\\\vdots&\vdots\\
0&0\\0&0
\end{matrix}
&
\begin{matrix}
 \dots \\ \dots\\ \phantom{\vdots}\\
  \dots \\ \dots
\end{matrix}&
\begin{matrix}
*&*\\ *&*\\\vdots&\vdots
\\0&\!H_{r_t}\!\\0&0
\end{matrix}
\end{array}\right],
\]
where each $H_{r_i}$ is a
nonsingular $r_i$-by-$r_i$
matrix. The matrix $A$ takes the
form
\begin{equation}\label{el.32qz}
\left[\begin{array}{cc|c|c}
\boxed{ \
\begin{matrix}
M_1&M_2&\dots&M_{t}&M_{t+1}
\\
0&0_{r_2}&\dots &0&0\\
\vdots&\vdots&\ddots&
\vdots&\vdots\\
 0&0&\dots &0_{r_{t}}&0\\
     0_{r_1}&0&\dots &0&0_{r_{t+1}}
\end{matrix}\;}&\!\!
\begin{matrix}
* &*\\0&H_{r_2}\\\vdots&\vdots\\0&0\\
0&0
\end{matrix}&
\begin{matrix}
 \dots \\ \dots\\ \phantom{\vdots}\\
  \dots\\ \dots
\end{matrix}&
\begin{matrix}
* &*\\ * &*\\\vdots&\vdots\\
0&H_{r_t}\!\!\\ 0&0
\end{matrix}
\end{array}\right].
\end{equation}
Replacing
$H_{r_{2}},\dots,H_{r_{t}}$ by
the identity matrices of the
same sizes and all $*$'s by the
zero matrices, we obtain
\eqref{el.33z} since
\eqref{el.32qz} reduces to
\eqref{el.33z} by those
transformations (a)--(c) from
Lemma \ref{s-res1} that preserve
$r_2,\dots,r_{t+1}$,
$M_1,\dots,M_{t+1}$.
\end{proof}

\subsection{Regularization algorithm
for similarity}\label{sub-r3z}

For any bangle
\begin{equation*}
A=\bigl[\,A_1\,|\dots|\,
A_{k-1}\ \boxed{\ A_k\ }\
A_{k+1}\,|\dots|\,A_t\,\bigr]
\end{equation*}
over $\mathbb F$, its
regularizing decomposition for
similarity can be constructed as
follows.

\begin{itemize}
  \item
First we apply subsequently the
left-hand reduction for
similarity to $A$ until obtain
\begin{equation*}
L(L\dots(L(A))\dots)
=\bigl[\,0_{m0}\,|\dots|\,
0_{m0}\ \boxed{\ B_k\ }\
B_{k+1}\,|\dots|\,B_t\,\bigr],
\end{equation*}
in which the first $k-1$ strips
have no columns.

  \item
Then we apply subsequently the
right-hand reduction for
similarity to
\begin{equation*}%\label{}
B=\bigl[\: \boxed{\ B_k\ }\
B_{k+1}\,|\dots|\,B_t\,\bigr]
\end{equation*}
until obtain
\begin{equation}\label{el.62xz}
R_s(R_s\dots(R_s(B))\dots)
=\bigl[\:\boxed{\ K\ }\
0_{n0}\,|\dots|\, 0_{n0}\,\bigr]
\end{equation}
with a nonsingular $K$.
\end{itemize}

Producing this reduction, we in
each step have deleted the
reduced parts of $A$. Let us
repeat the reduction preserving
all the reduced parts of $A$ and
denote the obtained bangle by
$\check A$. Clearly, $\check{A}$
is similar to $A$. Due to the
next theorem, $\check A$ is a
regularizing decomposition of
$A$ up to admissible
permutations of rows and
columns.

\begin{theorem}\label{tel.1z}
If $A$ is a bangle over a field
or skew field \/$\mathbb F$,
then $\check A$ reduces by
admissible permutations of rows
and columns to a regularizing
decomposition of $A$ for
similarity.
\end{theorem}

\begin{proof}
We give a constructive proof of
this theorem. By admissible
permutations of rows and
columns, $\check A$ is reduced
to a block-direct sum of the
bangle \eqref{e10.3} with $K$
from \eqref{el.62xz} and a
bangle $D$ in which each row and
each column contains at most one
$1$ and the other entries are
zero. Replacing $D$ in this
block-direct sum by $\Sigma_D$
from \eqref{43}, we obtain a
regularizing decomposition of
$A$ for similarity.
\end{proof}

\section{Proofs of Theorems
\ref{l10.1}, \ref{t0.01}, and
\ref{t.m}}\label{s-pr}

\begin{proof}[Proof of Theorem
\ref{l10.1}] (a) Let us prove
the statement (a) for
*congruence; its proof for
similarity is analogous.

Let $A$ be a bangle over
$\mathbb F$. In view of Theorem
\ref{tel.1}, $A$ possesses a
regularizing decomposition for
{\rm*\!}congruence, which is
obtained from $\hat A$ by
admissible permutations of rows
and columns.

Let $\Sigma_1$ and $\Sigma_2$ be
two regularizing decompositions
of $A$. Then $\Sigma_1
\overset{*}{\sim} \Sigma_2$. We
need to prove that
\begin{equation}\label{el.66}
\Sigma_1^{\mathrm{reg}}
\overset{*}{\sim}
\Sigma_2^{\mathrm{reg}} \quad
 \mathrm{and}\quad
\Sigma_1^{\mathrm{sing}}
 \overset{p}{\sim}
\Sigma_2^{\mathrm{sing}},
\end{equation}
where $\Sigma_i^{\mathrm{reg}}$
and $\Sigma_i^{\mathrm{sing}}$
are the regular and the singular
parts of $\Sigma_i$ ($i=1,2$).

If
\[
L(\Sigma_1) =R(\Sigma_1)
=\Sigma_1,
\]
then $\Sigma_1=
\Sigma_1^{\mathrm{reg}}$ and
\eqref{el.66} holds.

Let $L(\Sigma_1)\ne\Sigma_1$ or
$R(\Sigma_1) \ne\Sigma_1$.
Suppose for definiteness that
\begin{equation}\label{77}
L(\Sigma_1)\ne\Sigma_1.
\end{equation}

Each row and each column of
$\Sigma_i^{\mathrm{sing}}$
$(i=1,2)$ contains at most one
$1$, the other entries are zero.
Due to this property, the
reduction of $\Sigma_i$ to
\begin{equation}\label{kjk}
\Omega_i :={\cal
L}_k(L(\Sigma_i)).
\end{equation}
 of the form
\eqref{el.2} can be realized by
admissible permutations:
\begin{equation}\label{el.68}
\Sigma_i\overset{p}{\sim}\Omega_i;
\end{equation}
moreover, $L(\Sigma_i)$ is a
block-direct sum of a bangle of
the form \eqref{e10.3} and a
bangle, in which each row and
each column contains at most one
$1$, the other entries are zero.
By \eqref{43}, $L(\Sigma_i)$
reduces by admissible
permutations of rows and columns
to its regularizing
decomposition, so we may take
$\Omega_i$ such that
$L(\Sigma_i)$ is a regularizing
decomposition.

Since $\Sigma_1
\overset{*}{\sim} \Sigma_2$, we
have $\Omega_1 \overset{*}{\sim}
\Omega_2$, and so by
\eqref{el.5} and \eqref{kjk}
\begin{equation*}\label{lll}
L(\Sigma_1) \overset{*}{\sim}
L(\Sigma_2).
\end{equation*}
Due to \eqref{77}, the size of
$L(\Sigma_1)$ is less than the
size of $\Sigma_1$, reasoning by
induction we may assume that
\eqref{el.66} holds for
$L(\Sigma_i)$; that is,
\[
L(\Sigma_1)^{\mathrm{reg}}
\overset{*}{\sim}
L(\Sigma_2)^{\mathrm{reg}} \quad
 \mathrm{and}\quad
L(\Sigma_1)^{\mathrm{sing}}
 \overset{p}{\sim}
L(\Sigma_2)^{\mathrm{sing}}.
\]
Then
\[
\Omega_1^{\mathrm{reg}}
\overset{*}{\sim}
\Omega_2^{\mathrm{reg}} \quad
 \mathrm{and}\quad
\Omega_1^{\mathrm{sing}}
 \overset{p}{\sim}
\Omega_2^{\mathrm{sing}}
\]
since $\Omega_1$ and $\Omega_2$
have the form \eqref{el.2}. This
proves \eqref{el.66} due to
\eqref{el.68}.
 \medskip

(b) This statement follows from
(a) and Theorem \ref{bilin}.

(c) This statement follows from
(a) and the uniqueness of the
Jordan Canonical Form.
\end{proof}

\begin{proof}[Proof of Theorem
\ref{t0.01}] Let ${\cal G}$ be
one of sesquilinear forms
\begin{equation*}%\label{jsjb}
U\times V\rightarrow {\mathbb
F}, \qquad (V/U)\times
V\rightarrow {\mathbb F}.
\end{equation*}
Let us prove that the canonical
form problem for its matrix
$[A|B]$ (defined in
\eqref{0.22}) is the canonical
form problem under {\rm
*\!}congruence for the bangle
\begin{equation*}
\bigl[\,\boxed{\; A\; }\
B\,\bigr] \quad\text{or}\quad
\bigl[B\ \boxed{\;A\;}\;\bigr],
\end{equation*}
respectively, and so Theorem
\ref{t0.01} follows from Theorem
\ref{l10.1}.

 It suffices to
prove that a change of the basis
of $V$ reduces $[A|B]$ by
transformations
\begin{equation}\label{a.1}
[A\ B]\mapsto
  \begin{cases}
    S^*[A\ B]
\begin{bmatrix}
S&P\\0&Q
\end{bmatrix} &
 \text{if ${\cal G}\colon U\times V\to
{\mathbb F}$,}  \\[6mm]
     S^*[A\ B]
\begin{bmatrix}
S&0\\P&Q
\end{bmatrix}
&
 \text{if ${\cal G}\colon
(V/U)\times V\to {\mathbb F}$},
  \end{cases}
\end{equation}
in which $S$ and $Q$ are
nonsingular matrices and $P$ is
arbitrary.
\medskip

\noindent \emph{Case 1:} $[A|B]$
is the matrix of
\[
{\cal G}\colon U\times V\to
{\mathbb F},\qquad U\subset V,
\]
in a basis $e_1,\dots,e_n$ of
$V$ satisfying \eqref{0.21}. If
\begin{equation}\label{a.10}
f_{j}=e_1\rho_{1j}+\dots+
e_n\rho_{nj},\qquad j=1,\dots,n,
\end{equation}
is another basis of $V$ such
that $f_1,\dots,f_m$ is a basis
of $U$, then the change matrix
from $e_1,\dots,e_n$ to
$f_1,\dots,f_n$ has the form
\begin{equation*}\label{7.5}
R=[\rho_{ij}]=\begin{bmatrix}
S&P\\0&Q
\end{bmatrix},
\end{equation*}
where $S$ is the change matrix
from $e_1,\dots,e_m$ to
$f_1,\dots,f_m$ in $U$. Due to
\eqref{0.5}, the matrix $[A|B]$
reduces by transformations
\eqref{a.1}.
\medskip

\noindent \emph{Case 2:} $[A|B]$
is the matrix of
\[
{\cal G}\colon (V/U)\times V\to
{\mathbb F},\qquad U\subset V,
\]
in a basis $e_1,\dots,e_n$ of
$V$ satisfying \eqref{0.21}. If
\eqref{a.10} is another basis of
$V$ such that
$f_{m+1},\dots,f_n$ is a basis
of $U$, then the change matrix
from $e_1,\dots,e_n$ to
$f_1,\dots,f_n$ has the form
\begin{equation*}\label{7.5a}
R=[\rho_{ij}]=\begin{bmatrix}
S&0\\P&Q
\end{bmatrix},
\end{equation*}
where $S$ is the change matrix
from $e_1+U,\dots,e_m+U$ to
$f_1+U,\dots,f_m+U$ in $V/U$.
Hence, the matrix $[A|B]$
reduces by transformations
\eqref{a.1}.
\end{proof}

\begin{proof}[Proof of Theorem
\ref{t.m}] Let ${\cal A}$ be one
of linear mappings
\begin{equation*}%\label{jsj}
U\rightarrow V,\qquad
V\rightarrow U,\qquad
V/U\rightarrow V,\qquad
V\rightarrow V/U.
\end{equation*}
Let us prove that the canonical
form problem for its matrix
\begin{equation*}%\label{m.3a}
A_e= \begin{cases}
    \left[ \begin{array}{cc}
  A \\ \hline B
 \end{array}\right] & \text{if
$U\to V$ or $V/U\to V$,} \\[5mm]
   [A| B]
 & \text{if $V\to U$ or
 $V\to V/U$,}
  \end{cases}
\end{equation*}
(see \eqref{m.3}) is the
canonical form problem under
similarity for the bangle
\begin{equation*}
\Bigl[B^T\
\boxed{\;A^T\;}\;\Bigr],
    \quad
\bigl[\,\boxed{\; A\; }\
B\,\bigr],
   \quad
\Bigl[\,\boxed{\; A^T\; }\
B^T\,\Bigr],
   \quad\text{or}\quad
\bigl[B\ \boxed{\;A\;}\;\bigr],
 \end{equation*}
respectively, and so Theorem
\ref{t.m} follows from Theorem
\ref{l10.1}.

It suffices to prove that a
change of the basis of $V$
reduces $A_e$ by transformations
\begin{align}\label{qqqq1}
\begin{bmatrix}
  A \\ B
 \end{bmatrix}&
   \longmapsto
    \begin{bmatrix}
S^{-1}&*\\0&Q^{-1}
\end{bmatrix}
\begin{bmatrix}
  A \\ B
 \end{bmatrix}
S\quad
 \text{if ${\cal A}\colon U\to V$,}
     %%%%%%%%%%
          \\[2mm] \label{qqqq2}
[A\ B]&\longmapsto
    S^{-1}[A\ B]
\begin{bmatrix}
S& *\\0&Q
\end{bmatrix} \quad
 \text{if ${\cal A}\colon V\to U$,}
      %%%%%%%%
   \\[2mm]
      %%%%%%%%
      \label{qqqq3}
[A\ B]&\longmapsto
     S^{-1}[A\ B]
\begin{bmatrix}
S&0\\ *&Q
\end{bmatrix}
\quad
 \text{if ${\cal A}\colon
V\to V/U$},
      %%%%%%%%%%%
   \\[2mm] \label{qqqq4}
      %%%%%%%%%%%
\begin{bmatrix}
  A \\ B
 \end{bmatrix}&
   \longmapsto
   \begin{bmatrix}
S^{-1}&0\\ *&Q^{-1}
\end{bmatrix}
\begin{bmatrix}
  A \\ B
 \end{bmatrix}
S\quad
 \text{if ${\cal A}\colon
V/U\to V$},
\end{align}
in which $S$ and $Q$ are
nonsingular matrices and the
$*$'s denote arbitrary matrices.
\medskip

\noindent \emph{Case 1:} $A_e$
is the matrix of
\[
{\cal A}\colon U\to V
\quad\text{or}\quad {\cal
A}\colon V\to U,\qquad U\subset
V,
\]
in a basis $e_1,\dots,e_n$ of
$V$ satisfying \eqref{m.2}. If
\begin{equation}\label{ma.10}
f_{j}=e_1\rho_{1j}+\dots+
e_n\rho_{nj},\qquad j=1,\dots,n,
\end{equation}
is another basis of $V$ such
that $f_1,\dots,f_m$ is a basis
of $U$, then the change matrix
from $e_1,\dots,e_n$ to
$f_1,\dots,f_n$ has the form
\begin{equation*}
R=[\rho_{ij}]=\begin{bmatrix}
S&P\\0&Q
\end{bmatrix},
\end{equation*}
where $S$ is the change matrix
from $e_1,\dots,e_m$ to
$f_1,\dots,f_m$ in $U$. So the
matrix $A_e$ reduces by
transformations \eqref{qqqq1} or
\eqref{qqqq2}.
\medskip

\noindent \emph{Case 2:} $A_e$
is the matrix of
\[
{\cal A}\colon V/U\rightarrow
V\quad\text{or}\quad {\cal
A}\colon V\to V/U,\qquad
U\subset V,
\]
in a basis $e_1,\dots,e_n$ of
$V$ satisfying \eqref{m.2}.
 If
\eqref{ma.10} is another basis
of $V$ such that
$f_{m+1},\dots,f_n$ is a basis
of $U$, then the change matrix
from $e_1,\dots,e_n$ to
$f_1,\dots,f_n$ has the form
\begin{equation*}
R=[\rho_{ij}]=\begin{bmatrix}
S&0\\P&Q
\end{bmatrix},
\end{equation*}
where $S$ is the change matrix
from $e_1+U,\dots,e_m+U$ to
$f_1+U,\dots,f_m+U$ in $V/U$.
Hence, the matrix $A_e$ reduces
by transformations \eqref{qqqq3}
or \eqref{qqqq4}.
\end{proof}


\begin{thebibliography}{9}
%\end{thebibliography}

\bibitem{hor-ser_transpose}
R. A. Horn, V. V. Sergeichuk,
Congruences of a square matrix
and its transpose. \emph{Linear
Algebra Appl.} 389 (2004)
347--353.

\bibitem{hor-ser1_singul}
R. A. Horn, V. V. Sergeichuk, A
regularization algorithm for
matrices of bilinear and
sesquilinear forms, \emph{Linear
Algebra Appl.} 412 (2006)
380--395.

\bibitem{hor-ser1_canon}
R. A. Horn, V. V. Sergeichuk,
Canonical forms for complex
matrix congruence and
*congruence, \emph{Linear
Algebra Appl.} 416
(2006) 1010--1032.

\bibitem{naz}
L. A. Nazarova, A. V. Roiter, V.
V. Sergeichuk, V. M. Bondarenko,
Application of modules over a
dyad for the classification of
finite p-groups possessing an
abelian subgroup of index p and
of pairs of mutually
annihilating operators. \emph{J.
Soviet Math.} 3 (no. 5) (1975)
636--654.

\bibitem{ser_izvestiya}
V. V. Sergeichuk, Classification
problems for systems of forms
and linear mappings, {\it Math.
USSR Izvestiya,} 31 (no. 3)
(1988) 481--501.


\bibitem{ser}
V. V. Sergeichuk, Computation of
canonical matrices for chains
and cycles of linear mappings,
\emph{Linear Algebra Appl.} 376
(2004) 235--263.

\bibitem{doo}
P. Van Dooren, The computation
of Kronecker's canonical form of
a singular pencil, {\it Linear
Algebra Appl.}  27 (1979)
103--140.

\end{thebibliography}
\end{document}